\documentclass[12pt]{article}

\usepackage{latexsym,amsmath,graphics}
\usepackage[psamsfonts]{amssymb}

\newtheorem{theorem}{Theorem}
\newtheorem{theoremsansnumero}{Theorem}

\newtheorem{lemma}[theorem]{Lemma}
\newtheorem{proposition}[theorem]{Proposition}
\newtheorem{corollary}[theorem]{Corollary}

\def\mathF{\mathbb{F}}
\def\mathZ{\mathbb{Z}}
\def\mathN{\mathbb{N}}
\def\mathT{\mathbb{T}}
\def\mathQ{\mathbb{Q}}
\def\mathR{\mathbb{R}}
\def\mathM{\mathbb{M}}

\def\CalH{\mathcal{H}}
\def\CalT{\mathcal{T}}
\def\CalR{\mathcal{R}}
\def\CalS{\mathcal{S}}
\def\CalN{\mathcal{N}}
\def\CalU{\mathcal{U}}
\def\CalG{\mathcal{G}}
\def\CalF{\mathcal{F}}

\def\Gzero{\tilde{\Gamma}}
\def\Mzero{\tilde{M}}
\def\Fbig{\tilde{F}}
\def\deltabig{\tilde{\delta}}
\def\phibig{\tilde{\phi}}
\def\Amen{H}  

\newcommand\note[1]{}{\unskip}

\begin{document}

\thispagestyle{empty}
\title{AN UNCOUNTABLE FAMILY OF NON
ORBIT EQUIVALENT ACTIONS OF $\mathF_n$}
\author{Damien Gaboriau\thanks{C.N.R.S.} \and Sorin
Popa\thanks{
Supported in part by a NSF Grant
0100883}}
\date{}

\maketitle

\vsize=6.5in

\begin{abstract}  For each $ 2 \leq n \leq \infty$, we
construct an
uncountable family of free ergodic measure preserving actions
$\alpha_t$ of the free group $\mathF _n$ on the standard
probability space $(X, \mu)$ such that any two are non orbit
equivalent (in fact, not even stably orbit equivalent). These actions
are all ``rigid'' (in the sense of [Po01]), with the II$_1$
factors $L^\infty(X, \mu)\rtimes_{\alpha_t} \mathF _n$ mutually
non-isomorphic (even non-stably isomorphic) and in the class
$\CalH \CalT _{_{s}}.$
\end{abstract}

\begin{center} $\S $0.  \end{center}

Recall that two ergodic probability measure preserving
(p.m.p.) actions
$\sigma_{i}$ for $i=1,2$ of two countable groups $\Gamma_{i}$ on
probability measure
standard Borel spaces
$(X_{i},\mu_{i})$ are {\it orbit equivalent} (OE) if they define
partitions of the spaces into orbits that are isomorphic; more
precisely
if
there exists a measurable, almost everywhere defined isomorphism
$f:X_{1}\to X_{2}$
such that
$f_{*}\mu_{1}=\mu_{2}$ and the $\Gamma_{1}$-orbit of
$\mu_{1}$-almost
every $x\in
X_{1}$ is sent by $f$ onto the $\Gamma_{2}$-orbit of $f(x)$.

The theory of orbit equivalence, although underlying  the ``group
measure space
construction'' of Murray and von Neumann [MvN36], was born
with the work
of H.~Dye who proved, for example, the following striking result
[Dy59]: {\it Any two ergodic
p.m.p. free actions of $\Gamma_{1}\simeq \mathZ$ and
$\Gamma_{2}\simeq\oplus_{j\in \mathN } \mathZ /2\mathZ $ are
orbit equivalent.}

Through a series of works, the class of groups $\Gamma_{2}$
satisfying
Dye's theorem gradually increased until it achieved the necessary
and sufficient condition:
$\Gamma_{2}$ is infinite amenable [OW80]. In particular, all
infinite amenable
groups produce one and only one
ergodic p.m.p. free action up to orbit equivalence (see also
[CFW81] for a more general statement).

By using the notion of {\it strong ergodicity}, A.~Connes and
B.~Weiss proved that any non amenable group without Kazhdan
property (T) admits at least two non OE p.m.p. free ergodic
actions [CW80]. The first examples of groups with uncountably many
non OE free ergodic actions was obtained in [BG81], using a
somewhat circumstancial construction, based on prior work in
\cite{McD69}. Within the circle of ideas brought up by Zimmer's cocycle super-rigidity [Zi84],
certain Kazhdan property (T) lattices of Lie groups such as $\mathrm{SL}(n, \mathZ), n \geq 3$, were shown to admit uncountably many non OE free ergodic
actions as well (see [GG88]).
 It is only recently that we learned that this property
was shared by all infinite groups with Kazhdan property (T)
[Hj02-b], and thanks to another reason (bounded cohomology) by
many torsion free finitely generated direct products
$\Gamma_{1}\times\Gamma_{2}$, including non trivial ($l\geq 2$)
products of free groups like $\mathF_{p_1}\times\mathF
_{p_2}\times\cdots\times \mathF _{p_l}$, $p_{i}\geq 2$ [MS02].

On the other hand, the situation for the free groups themselves or
$\mathrm{SL}(2,\mathZ )$ remained unclear and arouse the interest of
producing more non OE free ergodic actions of $\mathF _n$ than
just the two given by [CW80] (more precisely in producing ways to
distinguish them from the OE point of view), reaching  the number
of 3 in 2001 [Po01], then 4 in 2002 [Hj02-b], for each finite $ n
\geq 2$.

Meanwhile, a certain rigidity phenomenon has been detected, behind
the notion of {\it cost} [Ga00] or $l^2$ Betti numbers [Ga01],
showing in particular that the rank of the free groups is an
invariant of OE, as well as the virtual Euler characteristic of
groups. As a consequence, all actions of finitely generated free
groups were shown to have trivial fundamental group ([Ga00] or
[Ga01, Cor. 5.7]).

But these techniques have no hold on various actions of the same
group. A new notion of rigidity for actions of a given group has
been found in [Po01] (Def. 5.10.1), called {\it relative property}
(T). It is an OE invariant, is satisfied by certain actions of the
free group $\mathF _n, 2 \leq n < \infty$, and not by others, and
led to the already mentioned third example by way of countability
of its outer automorphism group (5.3.3 or 8.7 in [Po01]). This
property will be crucial in our results and notice right away that
it is satisfied by the natural actions of $\mathrm{SL}(2, \mathZ )$ on the
$2$-torus $\mathT^2$ (cf. [Po01]). The main goal of this paper is
to prove the following:

\begin{theoremsansnumero} For each $2 \leq n \leq \infty$,
there exist uncountably many non orbit equivalent actions of the
free group $\mathF _n$, each being free ergodic and probability
measure preserving on a standard Borel space.
\end{theoremsansnumero}
Our result also applies to every virtually free group, and more
generally to every group that is virtually a non trivial free
product of infinite amenable groups. All actions in our examples
have the relative property (T), thus having at most countable
outer automorphism groups and countable fundamental group (trivial
in the case $\mathF_n$). Moreover, the uncountable family of
actions can be taken so that any two of them are non stably
equivalent. Similar results are obtained for the associated II$_1$
factors. All along the paper, we consider only standard Borel
spaces.

\begin{center} $\S $1.  \end{center}

  From the very beginning, the orbit equivalence theory has
been strongly
connected
with operator algebras through the ``group measure space
construction''
[MvN36, pp. 192-209], which was the first general way of
producing various
infinite dimensional von
Neumann algebras (even before the ``group von Neumann algebra''
[MvN43, Sect. 5.3]):
{\it each ergodic p.m.p. free action ${\alpha}$ of an infinite
countable
group $\Gamma$ on a standard Borel space
$(X,\mu)$ defines a type} II$_1$ {\it von Neumann  factor $M$
together with
the abelian subalgebra $A=L^{\infty}(X,\mu)\subset
L^\infty(X, \mu) \rtimes_{\alpha} \Gamma=M$ sitting in it as a
{Cartan subalgebra}. Moreover two ergodic actions define
isomorphic
inclusions $(A_{1}\subset M_{1})\simeq (A_{2}\subset M_{2})$
iff they are orbit equivalent} [FM77]. This fact establishes
a useful equivalence between the study of actions of groups,
up to orbit
equivalence, and the study of the inclusions
of Cartan subalgebras into II$_1$ factors entailed by such
actions,
up to isomorphism of inclusions.

The notion of relative property (T) for an action of a group in
[Po01] is in fact expressed in terms of its associated Cartan
subalgebra, and the terminology used in this functional analytical
framework is {\it relatively rigid Cartan subalgebra} (or {\it
rigid inclusion}) $A \subset M$. Inspired by the notion of
relative property (T) of Kazhdan-Margulis for inclusions of groups
[Ma82] (see also [dHV89, p. 16]) and by the Connes-Jones
definition of property (T) for II$_1$ factors [CJ85],  this
property requires that whenever a completely positive map on $M$
(the operator algebra analogue of a positive definite function) is
close to the identity on some ``critical'' finite subset of $M$,
then it must be uniformly close to the identity on $A$.

If a Cartan subalgebra $A \subset M$ comes from a group von
Neumann algebra construction corresponding to an inclusion of
groups $H \subset G=\Gamma \ltimes H$ then $A \subset M$ is rigid
in the sense of Definition 4.2 in [Po01] iff $H \subset G$ has the
relative property (T) (cf. 5.1 in [Po01]). Since $\mathZ ^2
\subset \mathrm{SL}(2, \mathZ)  \ltimes \mathZ ^2$ has the relative
property (T) (cf. [Ma82]), this observation gives rise to an
important example of an action with the relative property (T),
namely the standard action $(\mathrm{SL}(2, \mathZ ),\mathT^2)$. More
generally, M. Burger has shown that for any non-amenable subgroup
of automorphisms of $\mathZ^2$, $\Gamma \subset \mathrm{SL}(2, \mathZ )$,
the induced inclusion $\mathZ ^2 \subset \Gamma \ltimes \mathZ ^2$
has the relative property (T). So in fact by (5.1 in [Po01]) all
actions $(\Gamma, \mathT ^2)$ have the relative property (T).
Equivalently, $A \subset M$ is a rigid inclusion of von Neumann
algebras, where $M$ is the group von Neumann factor of the natural
semi-direct product $\Gamma \ltimes \mathZ ^2$ and $A$ is the
Cartan subalgebra $A=L^{\infty}(\mathT ^{2})$, identified via
Fourier transform with the von Neumann algebra of the group
$\mathZ ^{2}$.

For the proof of the Theorem, we start with the above example of
rigid inclusion $A \subset M$ in which $\Gamma=\mathF _n$. Then we
use Dye's theorem to write the equivalence relation $\CalR _1$
implemented by the first of the generators of $\mathF _n$ (which
one may assume acts ergodically) as a limit of an increasing one
parameter family of ergodic hyperfinite equivalence relations
$\CalR _{1,t}$, for  $t\in (0,1]$. Together with $\mathF _{n-1}$,
each of the $\CalR _{1,t}$ gives rise to a free ergodic action
$\alpha_t$ of $\mathF _n$ on $\mathT ^2$, with strictly increasing
group measure space factors $A \subset M_t \subset M$ and $M_t
\nearrow M$. By a rigidity property [Po01, 4.8], there is a value
$c$ of the parameters after which $(c<t\leq 1)$ the actions
$\alpha_{t}$ have relative property (T), equivalently, $A \subset
M_t$ are rigid.

We then use the rigidity of each of these $A \subset M_t$ to prove
that if an endomorphism of $M_t$ is close enough to the identity,
then it has to be onto. This is done by analogy with the proof
that Out$(A \subset M_t)$ is countable in [Po01, 4.7], along the
lines of Connes' pioneering rigidity result on the countability of
Out$(N)$ for II$_1$ factors $N$ with the property (T) [Co80].
Thus, by a separability argument as in [Co80], [Po86], if
uncountably many inclusions $(A \subset M_t)$ are isomorphic then
one gets an endomorphism $M_{s}\simeq M_{t}\subset M_{s}$ close to
the identity, thus onto, contradicting $M_t \neq M_s$. The
partition of the interval of parameters $[c,1]$ into the subsets
where the factors $M_{t}$ are isomorphic is then composed of
finite or countable subsets. We mention that in fact when proving
the statement involving stable orbit equivalence, the above
argument becomes  slightly more complicated, with ``isomorphism''
(resp. ``endomorphism'') replaced by surjective (resp.
non-surjective) isomorphism between amplifications of suitable
algebras.

\begin{center} $\S $2.  \end{center}

Set $I=(0,1)\cap \mathQ $ and $I_{t}=I\cap (0,t]$, for $t\in
(0,1]$;
set  $G=\oplus_{i\in I}\mathZ /2 \mathZ $ and
$G_{t}=\oplus_{i\in I_{t}}
\mathZ  /2\mathZ $, for $0<t\leq 1$. This gives a strictly
increasing continuum
of subgroups of $G$.

\begin{proposition}\label{prop:1}
Let ${\alpha}$ be a free p.m.p. action of \ ${\Gamma_{1}} *
{\Gamma_{2}}$, with
${\Gamma_{1}}$ an infinite amenable group acting ergodically. Then

1$^{\circ}.$
The action ${\alpha}$ is orbit equivalent to a free ergodic
action $\sigma_{1}$
of the group $G*{\Gamma_{2}}$.

2$^{\circ}.$
The restrictions $\sigma_{t}$
of $\sigma_{1}$ to the subgroups $G_{t}*{\Gamma_{2}}$ define a
strictly increasing family of ergodic equivalence relations
${\CalS }_{t}$ indexed by $t\in (0,1]$, such that for any $s$
and for
any increasing
sequence $t_n \nearrow s$ one has
${\CalS }_{s}=\cup_{n}{\CalS }_{t_n}$.

3$^{\circ}.$ The action $\sigma_{t}$ of $G_{t}*{\Gamma_{2}}$
is OE to
some ergodic
free action $\alpha_{t}$ of ${\Gamma_{1}} * {\Gamma_{2}}$.

4$^{\circ}.$ The associated Cartan subalgebra inclusions define
a strictly increasing
family of subfactors $A \subset M_t \subset M_1$,
$0 < t \leq 1$, such that $M_t  =
A \rtimes_{\alpha_t} (\Gamma_1 * \Gamma_2)$,
$M_1= A \rtimes_{\alpha} (\Gamma_1 * \Gamma_2)$,
and $\overline{\cup_{t<s} M_t}^w= M_s$,  for each $s\in (0,1]$.
\end{proposition}

\begin{lemma}\label{lem:2}
Let $\Gamma=B*D$ and $\Lambda=C*D$
be free products such that $B$ and $C$ are infinite amenable
groups.
Every p.m.p. free action $\alpha$
of $\Gamma$, where $B$ acts ergodically,
is orbit equivalent to a free action $\sigma$
of $\Lambda$, where (up to a null set) $\sigma_{|C}$ is any
prescribed
free ergodic p.m.p. $C$-action and $\sigma_{|D}
\simeq\alpha_{|D}$.
\end{lemma}
\noindent {\it Proof}. By Ornstein-Weiss' [OW80] and Dye's
[Dy59] theorems,
the free action of $B$ on $(X_1,\mu_1)$ is OE to the
prescribed $C$-action
on $(X_2,\mu_2)$, as witnessed by $f:X_1 \to X_2$. Conjugating
$\alpha_{|D}$
by $f$ leads to an action of the free product $C*D$, which is
OE to $\alpha$.
   The equivalence relations defined by the actions,
${\CalS}_{1}={\CalS }_{D}$ and ${\CalS }_{2}={\CalS }_{B}={\CalS
}_{C}$,
are in free product in the sense of [Ga00, D\'ef. IV.9], so that a
non trivial reduced word in the free product $C*D$ has
obviously a set of
fixed point of measure zero: the $\Lambda$-action $\sigma$ is
free. \hfill
Q.E.D.

\medskip

\noindent
{\it Proof of the Proposition~\ref{prop:1}}.
By Lemma~\ref{lem:2}, ${\alpha}$ is orbit equivalent to a free
action
$\sigma_{1}$
of $G*{\Gamma_{2}}$, where the $G$-action is its Bernoulli
shift on
$\{0,1\}^{G}$.
The restriction of $\sigma_1$ to $G_{t}$ remains ergodic, and
the restriction
$\sigma_{t}$
of $\sigma_1$ to $G_{t}*{\Gamma_{2}}$ is free
ergodic. The properties of the family
${\CalS }_{t}$
are now obvious, and another application of Lemma~\ref{lem:2}
shows that each
$\sigma_{t}$ is OE to
a free action of ${\Gamma_{1}} * {\Gamma_{2}}$.
The operator algebra statement 4$^{\circ}$ follows
immediately,
due to the above and [FM77].
\hfill Q.E.D.

\bigskip
Recall that a {\it Cartan subalgebra} $A$ in a factor $M$
is a maximal abelian subalgebra whose normalizer $\CalN
_M(A)\overset{\mathrm{def}}{=}
\{u\in {\CalU (M)} \mid u A u ^{*}=A \}$ generates $M$ as von
Neumann algebra
(where ${\CalU (M)}$ is the group of unitary elements of $M$).
Also, recall from 4.2 in [Po01] that a Cartan subalgebra $A$
in a II$_1$
factor $M$ is called {\it relatively rigid} if for any
$\varepsilon > 0$ there exist a finite set $F\subset  M$ and
$\delta > 0$ such that if $\phi: M \rightarrow M$ is a completely
positive map with $\tau\circ \phi \leq \tau, \phi(1)\leq 1$ and
   $\|\phi(x) - x \|_2 \leq \delta, \forall x\in F$, then
$\|\phi(a)-a\|_2 \leq \varepsilon, \forall a\in A$, $\|a\| \leq
1$.

\vskip .05in In [4.7, Po01], it was proved that if $M$ is a type
II$_1$ factor with $A \subset M$ a relatively rigid Cartan
subalgebra then the subgroup of automorphism $\CalG _A \subset
{\text{\rm Aut}}(M)$ that are inner perturbations of automorphisms
leaving $A$ pointwise fixed is open and closed in Aut$(M)$, thus
giving rise to a countable quotient group Aut$(M)/\CalG _A$. In
particular, if $A \subset M$ comes from a free ergodic p.m.p.
action ${\alpha}$ of a group $\Gamma$ and $\CalR $ denotes the
corresponding equivalence relation implemented by ${\alpha}$ on
the probability space, then Out$(\CalR )$ is countable. The next
Proposition generalizes this result to the case of (not
necessarily unital) endomorphisms of $M$ and of its
amplifications.

\begin{lemma}\label{prop:3}
  Let $M$ be a type ${\text{\rm II}}_1$ factor
with $A \subset M$ a relatively rigid Cartan subalgebra. For each
$x \in M$ denote by $i(x)=x \oplus 0\in \mathM_{2 \times 2}(M)$
and $p_0 = i(1_M)=1\oplus 0$. There exist $F \subset M$
finite and $\delta
>  0$ such that if $\theta : i(M) \rightarrow p\mathM_{2
\times 2}(M)p$ is a unital $*$-isomorphism satisfying
$\|\theta(i(x))-i(x)\|_2 \leq \delta$, $\forall x\in F$, for some
projection $p \in \mathM_{2 \times 2}(M)$, then there exists
$p'\in \theta(i(A))'\cap p\mathM_{2 \times 2}(M)p$ and a partial
isometry $u \in  \mathM_{2 \times 2}(M)$ such that
\begin{center}
$u^*u=p_0, uu^*=p', \theta(i(M))p'=p'\mathM_{2 \times 2}(M)p',
\theta(i(a))p'=ui(a)u^*, \forall a\in A.$
\end{center}
If in addition $\theta(i(M))'\cap p\mathM_{2 \times 2}(M)p = \mathbb C p$, then $\theta$ is surjective and $\tau(p)=1/2$, i.e., $\tau$
preserves the trace.
\end{lemma}
\noindent {\it Proof}. Since $\theta$ is a unital $*$-isomorphism,
it is completely positive. Thus, if $\phi : M \rightarrow M$ is
defined by $\phi (x) = p_0 \theta (x)p_0, x\in M$, then $\phi$ is
completely positive and for all $x\in M$ we have:
\begin{eqnarray}
\|\phi(x) - x\|_2 & \leq & \|\theta(i(x)) - i(x)\|_2 +
\|\theta (i(x))
- p_0\theta (i(x)) p_0\|_2 \nonumber \\
& \leq & \|\theta(i(x)) - i(x)\|_2 + \|\theta (i(x)) -
\theta(i(x))p_0\|_2 \nonumber \\
& &
\hskip 100pt
+ \|\ \theta(i(x)) p_0 - p_0\theta (i(x))
p_0\|_2 \nonumber \\
&  \leq & \|\theta(i(x)) - i(x)\|_2 + 2 \| \theta(p_0)-
p_0\|_2.\nonumber
\end{eqnarray}

By the definition of relative rigidity it follows that there exist
$F\subset M$ finite, with $1_M\in F$, and $\delta > 0$, with
$\delta \leq 1/8$, such that if a subunital subtracial completely
positive map $\phi$ on $M$ satisfies $\|\phi(x)-x\|_2 \leq \delta,
\forall x\in F$, then $\|\phi(v)-v\|_2 \leq 1/4$ for all elements
$v$ in the unitary group $\CalU (A)$ of $A$. But then we also have
\begin{eqnarray}
\|\theta (i(v)) - i(v)\|_2 &\leq& \|\theta (i(v)) - \theta (i(v))
p_0\|_2 + \|\theta(i(v)) p_0 - p_0\theta (i(v)) p_0\|_2
\nonumber\\
& & \hskip 220pt
+ \|\phi(v) - v\|_2\nonumber\\
& \leq & 2\|\theta (p_0) - p_0\|_2 + \|\phi(v)-v \|_2\ \nonumber\\
&  \leq &
2 \cdot 1/8
+ 1/4\ =\ 1/2.\nonumber
\end{eqnarray}
Thus, if $h$ denotes the unique element of minimal Hilbert norm
$\|\cdot \|_2$ in the weakly compact convex set
$\overline{\text{\rm co}}^w \{\theta(i(v))i(v)^* \mid v \in \CalU
(A)\}\subset p\mathM_{2 \times 2}(M)p_0$ then $\|h - p_0\|_2 \leq
1/2$. Also, since $\|\theta(i(v))hi(v)^*\|_2 = \|h\|_2, \forall v
\in \CalU (A)$, by the uniqueness of $h$ it follows that
$\theta(i(v))h=hi(v), \forall v \in \CalU (A).$ By a standard
functional calculus trick (as in group representations), if
$u_0\in \mathM_{2 \times 2}(M)$ is the partial isometry in the
polar decomposition of $h$ then $u_0 \neq 0$,
$u_0i(a)=\theta(i(a))u_0, \forall a\in A$, and $u_0^*u_0\in
i(A)'\cap i(M)=i(A)$, $[u_0u_0^*, \theta(i(A))]=0$.

We can thus apply the above to $a=i^{-1}(u_0^*u_0)\in A$ to get
$\theta(u_0^*u_0)u_0=u_0 u_0^*u_0=u_0$. But this implies
$\theta(u_0^*u_0)u_0u_0^*=u_0u_0^*$, so that $\theta(u_0^*u_0)
\geq u_0u_0^*$. Note that in particular this implies that $\theta$
doesn't ``shrink'' the trace $\tau$, i.e., $\tau(\theta(x)) \geq
\tau(x), \forall x\in p_0\mathM_{2 \times 2}(M)p_0$.

Since the normalizer $\CalN _M(A)$ of $A$ in $M$ acts ergodically
on $A$, there exist partial isometries $v_1, v_2, \ldots, v_n$ in
$M$ such that $v_j v_j^*, v_j^* v_j \in A$, $v_jAv_j^* = Av_j
v_j^*$, $v_j^* v_j \leq u_0^* u_0$ and $\sum_j v_j v_j^* = 1$.
Then $u = \sum_j \theta(i(v_j))u_0i(v_j)^*$ is easily seen to be a
partial isometry satisfying $ui(a)= \theta(i(a))u, \forall a \in
A$, with $u^*u = p_0$, $p'=uu^* \in \theta(i(A))' \cap p\mathM_{2
\times 2}(M)p$. By spatiality it follows that $\theta(i(A))p'$ is
a Cartan subalgebra (in particular maximal abelian) in $p'\mathM_{2 \times 2}(M)p'$.

Let $B = \theta(i(A))'\cap p\mathM_{2 \times 2}(M)p$. Since
$\theta(i(w))$ normalizes $\theta(i(A))$ for  $w \in \CalN(A)$ it
follows that $\theta(i(w))$ normalizes $B$ as well. Also,
$\theta(i(A))p'$ maximal abelian implies $p'Bp'= \theta(i(A))p'$.
Thus, if $p_1 \in B$ denotes the maximal abelian central
projection of $B$ then $p_1 \leq p'$ and $\tau(p'-p_1) \leq
\tau(p-p')$. Since $\tau(p')=\tau(p_0)=1/2$ and
$$
\tau(p)-\tau(p_0) \leq \|p-p_0\|_1 \leq 2 \|p-p_0\|_2 \leq 2
\delta \leq 1/4,
$$
it follows that $1/2 - \tau(p_1) \leq 1/4$, so that $p_1 \neq 0$.
But since $B$ is normalized by $\theta(i(\CalN(A)))$, $p_1$ is
left invariant by Ad$(w')$ for $w' \in \theta(i(\CalN(A)))$. Thus,
$[p_1, \theta(i(M))]=0$, so in particular $a \mapsto ap_1$, $a \in
\theta(i(A))$ is an isomorphism. But $p'Bp'=\theta(i(A))p'\simeq
A$ as well. Thus, $p_1 = p'$.

In particular, this shows that if we denote
$\theta_0(x)=u^*\theta(x)u, x\in M$, then $\theta$ is an
isomorphism. Moreover, $\theta_0(a)=a, \forall a\in A$. Then we
get:
$$
waw^*= \theta_0(waw^*)=\theta_0(w)\theta_0(a)\theta_0(w^*) =
\theta(w)a\theta_0(w^*)
$$
for all $ a\in A$, implying that $\theta_0(w)^*w \in A'\cap M=A$.
Hence, $w\in \theta_0(w)A$. This shows that $\CalN _M(A) \subset
\theta_0(M)A=\theta_0(M)$, so that $M \subset \theta_0(M)$, thus
$\theta_0(M)=M$. But this implies $\theta(M)=M$ as well. \hfill
Q.E.D.

\bigskip

To state the next result, recall from \cite{MvN43} that if $M$ is
a II$_1$ factor and $s
>  0$ then the $s$-amplification of $M$ is the II$_1$ factor $M^s=
p\mathM_{n \times n}(M)p$, where $n$ is some integer $\geq s$ and
$p$ is a projection in $\mathM_{n \times n}(M)$ of (normalized)
trace equal to $s/n$. This factor is defined up to isomorphisms
implemented by unitary elements in large enough matrix algebras
over $M$. Note also that in case $M$ has a Cartan subalgebra $A
\subset M$, then $p$ can be taken in the Cartan subalgebra
${\oplus^n}A$ of $\mathM_{n \times n}(M)$ given by diagonal
elements with entries in $A$ and $A^t = p({\oplus^n}A)p$ is a
Cartan subalgebra of $M^t$.

Two factors of type II$_1$, $M_1, M_2$ are {\it stably isomorphic}
if there exist $s > 0$ such that $M_1$ is isomorphic to $M_2^{s}$.
We then write $M_1 \cong M_2$, while the symbol $\simeq$ will be
used to denote usual isomorphism. Recall that the {\it fundamental
group} of $M$ is the subgroup $\CalF(M):=\{s:M^s\simeq M\}$ of
$\mathR^{*}_{+}$.

\begin{proposition}\label{prop:4}
Let $M$ be a type $\text{\rm II}_1$ factor with a relatively rigid
Cartan subalgebra $A \subset M$. Then the fundamental group of $M$
is at most countable. If in addition  $A \subset M_t \subset M$,
$t\in (0,1]$, is a strictly increasing family of subfactors such
that $M = \overline{\cup_{t<1}M_t}$, then we have:

$1^\circ$. There
exists $c < 1$ such that for all
$c \leq t \leq 1$,  $A \subset M_t$ are relatively rigid
Cartan subalgebras.

$2^\circ$. If $c$ satisfies $1^\circ$, then for each $t \in [c,1]$
the class of stable isomorphism $\{s \mid M_s \cong M_t\}$ is at most countable.

Thus, the quotient set of parameters $[c,1]/\sim$ is uncountable.
\end{proposition}

\noindent {\it Proof}. Item $1^\circ$ is part of $4.5.4^\circ$ in
[Po01]. To prove $2^\circ$, assume  there exist $t_0 \in [c,1]$
and an uncountable subset $L \subset [c,1]$ such that $\theta_t :
M_{t_0} \rightarrow M_t^{s(t)}$, $\forall t \in L$, for some
$s(t)> 0$ and some surjective isomorphisms $\theta_t$. We may
further assume $n_0^{-1} \leq s(t) \leq n_0, \forall t \in L$, for
some integer $n_0\geq 1$.

Let $F\subset M_{t_0}, \delta > 0$ be so that if $\theta$ is a
unital  isomorphism of $i(M_{t_0})=M_{t_0}\oplus 0$ into $p\mathM_{2 \times 2}(M_{t_0})p$ satisfying $\theta(i(M_{t_0}))'\cap
p\mathM_{2 \times 2}(M_{t_0})p = \mathbb Cp$ and
$\|\theta(i(x))-i(x)\|_2 \leq \delta, \forall x\in F$, as in
Lemma~\ref{prop:3}. Then $\theta$ is onto and preserves the
trace (cf. Lemma~\ref{prop:3}). By isomorphism, it follows
that for each $t \in L$ we have: If $\theta$ is a unital
isomorphism of $i(M_t^{s(t)})$ into $p\mathM_{2 \times
2}(M_{t}^{s(t)})p$ with range having trivial relative commutant
and satisfying $\|\theta(i(x))-i(x)\|_2 \leq \delta, \forall x\in
\theta_t(F)$, then $\theta$ is surjective and preserves the trace.

We will now regard the isomorphisms $\theta_t, t\in L$, as
isomorphisms of $M_{t_0} \oplus 0 \subset M \oplus 0$ into $\mathM_{n_0 \times n_0}(M)$, in the obvious way, where $M \oplus 0$ is
the algebra of $n_0\times n_0$ matrices of the form $x \oplus 0$,
having $x\in M$ as upper left entry and $0$ elsewhere.

By the separability of the Hilbert space $L^2(\mathM_{n_0 \times
n_0}(M), \tau)$ and the uncountability assumption on $L$, it
follows that there exist $t_1\neq t_2 \in L$ such that
$\|\theta_{t_1}(x)-\theta_{t_2}(x)\|_2 \leq \delta/n^2_0, \forall
x\in F$. Assume $t_1 < t_2$. Denote $\theta=\theta_{t_1} \circ
\theta_{t_2}^{-1} : M^{s(t_2)}_{t_2} \simeq M_{t_1}^{s(t_1)}
\subset M_{t_2}^{s(t_1)}$ and note that for $\delta \leq 1/8$ we
may assume $\theta$ is an isomorphism of $i(M_{t_2}^{s(t_2)})$
into $\mathM_{2 \times 2}(M_{t_2}^{s(t_2)})$ satisfying
$\|\theta(i(x))-i(x)\|_2 \leq \delta,$ $\forall x\in
\theta_{t_2}(F) \subset \theta_{t_2}(M_{t_0})=M_{t_2}^{s(t_2)}$.
Note also that since $M_{t_1}'\cap M_{t_2} = \mathbb C$, the range of
$\theta$ has trivial relative commutant.

By the above observation (consequence of Lemma~\ref{prop:3})
this implies $\theta$ is onto and trace preserving. Thus
$s(t_1)=s(t_2)$ and $M_{t_1} = M_{t_2}$, a contradiction.

Finally, note that the above argument applied to just the factor
$M$ and its amplifications shows that $\CalF(M)$ is countable.
\hfill Q.E.D.

\vskip .05in For the next result recall that the {\it fundamental
group} $\CalF(\CalR )$ of an ergodic equivalence relation $\CalR $ is the
set of all $t > 0$ such that the $t$-amplification $\CalR^t$ of
$\CalR $ is orbit equivalent to $\CalR $. Also, two equivalence
relations $\CalR_1, \CalR_2$ are {\it stably orbit equivalent} ($\CalR_1
\cong \CalR_2$) if there exists $s > 0$ such that $\CalR_1,
\CalR_2^{s}$ are orbit equivalent. Note that if $M, M_1, M_2$ are
the II$_1$ factors associated to $\CalR, \CalR_1, \CalR_2$,
respectively, then $\CalF(M) \supset \CalF(\CalR)$ and $\CalR_1
\cong \CalR_2$ implies $M_1 \cong M_2$.

\begin{corollary}\label{cor:5} For each $ 2 \leq n \leq
\infty$ there exists an uncountable family of non stably orbit
equivalent free ergodic p.m.p. actions $\alpha_{t}$ of $\mathF _n$
such that moreover the corresponding equivalence relations $\CalR
_{\alpha_t, \mathF _n}$ have at most countable fundamental group
(trivial in the case $n < \infty$) and at most countable outer
automorphism group.
\end{corollary}

\smallskip
\noindent {\it Proof}. For each $ 2 \leq n \leq \infty$, let
$\mathF_{n}\subset \mathrm{SL}(2, \mathZ )$ be a non-cyclic free
subgroup, and
consider its standard action ${\alpha}$ on the $2$-torus
$\mathT^{2}$. The following lemma~\ref{lem:6} gives a base of
$\mathF
_{n}$ with at
least one ergodically acting element, thus allowing the use of
Proposition~\ref{prop:1}. The family of Cartan subalgebra
inclusions given by
Proposition~\ref{prop:1} satisfies Proposition~\ref{prop:4}, since
$\alpha$ has the
relative property (T), {\it i.e.}
   $A=L(\mathZ ^2)=L^\infty(\mathT ^2,\lambda)$ is a relatively
rigid Cartan subalgebra of $M_1=L(\mathZ ^2\rtimes \mathF _{n})=
L(\mathT ^2, \lambda) \rtimes_{{\alpha}} \mathF _{n}$ (by [Bu91]
and [5.1--5.2, Po01]). But non isomorphic factors necessarily come
from non OE actions.

Thus, the uncountable  family of non stable equivalent actions,
each one with a countable fundamental group, can be taken to be
the $\alpha_{t}$ of Proposition~\ref{prop:1}, for a set of
parameters representative of the classes $[c,1]/\sim$ in
Proposition~\ref{prop:4} (2$^{\circ}$).
The triviality of the fundamental group in the case $n < \infty$
is a consequence of \cite[Cor. 5.7]{Ga01}.
The countability of the outer automorphism group follows from [Cor. 4.7, Po01]. \hfill
Q.E.D.

\begin{lemma}\label{lem:6}
If $\Gamma$ is a non cyclic free subgroup of $\mathrm{SL}(2, \mathZ )$,
then it admits a base of elements whose standard actions on the
$2$-torus $\mathT ^2$ are individually ergodic.
\end{lemma}

\noindent {\it Proof}. The criterion for an element to act
ergodically is that it doesn't have a root of unity as eigenvalue,
or equivalently that it is hyperbolic when looking at its action
on the circle, boundary of the hyperbolic disk. Fact: if $a$ and
$b$ are two non commuting elements of infinite order, where $a$ is
parabolic, then there exists a power $k\in \mathZ $ such that
$b^{k}a$ is hyperbolic: let $P_{a}$ be the unique fixed point of
$a$, let $P_{b}\not=P_{a}$ the first fixed point of $b$ met on the
circle from the attracting side of $P_{a}$ and let $I$ be the
corresponding interval between them (the case of a common  unique
fixed point is excluded by non commutativity, then non
solvability, of the subgroup generated by $a$ and $b$). Let $J$ be
a closed interval strictly contained in $I$ with end point
$P_{b}$. Then $aJ$ is contained in the interior of $I$, and up to
taking negative powers $P_{b}$ is attracting so that a big enough
powers $b^{k}$ satisfies $(b^{k}a) J\subset$ Int$(J)$, which is
equivalent to hyperbolicity. Finally, starting from a base of
$\Gamma$, one can successively transform it accordingly to the
fact so as to produce a base satisfying the lemma. \hfill Q.E.D.

\vskip .05pt

For next statement, recall that a group is virtually of a certain kind if it contains a finite index subgroup of that kind. The non triviality
hypothesis aims to avoid $\Lambda\simeq \Lambda*\{1\}$.

\begin{corollary}\label{cor:7} Corollary~\ref{cor:5} applies
equally well with, instead of the group $\mathF _{n}$, any
countable discrete group $\Gamma$ that is virtually a  non trivial
free product $\Lambda$ of infinite amenable groups. More generally,
if $\Gamma=\Lambda* H$ admits a free p.m.p. action $\alpha$ with the relative property (T), where $H$ is infinite amenable and acts ergodically, then $\Gamma$ has an uncountable family of non stably orbit
equivalent free ergodic p.m.p. actions $\alpha_{t}$,
such that $\CalR_{\alpha_t, \Gamma}$ have at most countable fundamental group
(trivial in the case $\beta_1(\Lambda)$ is finite) and at most countable outer
automorphism group.
\end{corollary}
Note that in order to ensure that the above action $\alpha$ has the relative property (T), it is enough to get this property for some subgroup of $\Gamma$.
This is, of course, trivial from the definition of the relative property (T)
(\cite[4.6.2]{Po01}).

\medskip
\noindent {\it Proof}. The first part relies on the fact that actions of
$\Gamma$ may be produced from actions of $\mathF _n$. Let
$\Lambda$ be a finite index subgroup of $\Gamma$, then the usual
suspension construction associates, to each free p.m.p. action of
$\Lambda$ on some standard Borel space $X$, a stably orbit
equivalent free p.m.p. action of $\Gamma$ (given by right
multiplication on itself) on the diagonal action quotient
$\Lambda\backslash(X\times \Gamma)$. If the $\Lambda$-action is
ergodic, then so is the $\Gamma$-action. Two non OE
$\Lambda$-actions lead to two non OE $\Gamma$-actions.

By hypothesis, $\Lambda$ may be chosen to be a free product
$H_1*H_2*\ldots$ of $n\in \{2,3,\ldots, \infty\}$ infinite
amenable groups. Start with an uncountable collection of free
p.m.p. actions of $\mathF _n$ given by corollary \ref{cor:5}.
Observe that in that construction one may assume that $\mathF _n$
is given with a base of elements that act individually in an
ergodic manner.

Then up to removing at most countably many times a null set, as in
lemma~\ref{lem:2}, one can replace each cyclic factor in the
decomposition of $\mathF _n$ by one $H_i$. This leads to
uncountably many free p.m.p actions of $\Lambda$.\hfill

The last part is an immediate consequence of Propositions 1 and 4.
Again, the triviality of the fundamental group in the case
$\beta_1(\Lambda) < \infty$ follows from \cite[Cor. 5.7]{Ga01}.
Q.E.D.

\vskip .05in

We'll now derive some applications of the above results to type
II$_1$ factors. To this end, recall from Section 6 in [Po01] that
a Cartan subalgebra $A$ of a II$_1$ factor $M$ is called
HT$_{_{s}}$ if in addition to $A \subset M$ being a rigid
inclusion, one has that $M$ has the {\it property H relative to}
$A$, i.e., there exists a sequence of unital completely positive
$A-A$ bimodular maps $\phi_n$ on $M$ that tend to the identity in
the point $\|\cdot \|_2$-topology and are compact relative to $A$
({\it cf.} 2.1 in [Po01]). A II$_1$ factor $M$ is in {\it the
class} $\CalH \CalT _{_{s}}$ if it has such a HT$_{_{s}}$ Cartan
subalgebra.

The property H of $M$ relative to $A$  is automatically satisfied
whenever $M$ comes from a ``group measure space construction'',
$M=L^\infty(X, \mu) \rtimes_{\alpha} \Gamma$, with $\Gamma$ a
group satisfying Haagerup's compact approximation property and
${\alpha}$ a free ergodic measure preserving action of $\Gamma$ on
the probability space $(X, \mu)$. Thus, since all subgroups
$\Gamma$ of $\mathrm{SL}(2, \mathZ )$ do have Haagerup's property (cf.
[Ha79]) and $(\Gamma \ltimes \mathZ^2, \mathZ^2)$ has the relative
property (T) for all non-amenable $\Gamma \subset \mathrm{SL}(2, \mathZ)$
(cf. [Bu91]), it follows by Lemma 6 that if $\Gamma=\mathF_n$ for
some $2 \leq n \leq \infty$ then $A \subset M=L^\infty(\mathbb T^2,
\lambda) \rtimes \Gamma$ is a HT$_{_{s}}$ Cartan subalgebra of the
HT$_{_{s}}$ factor $M$. We then have:

\begin{corollary}\label{cor:8} For each $ 2 \leq n \leq
\infty$ there exists an uncountable family of mutually non stably
isomorphic II$_1$ factors $M_t$ of the form $L^\infty(X, \mu)
\rtimes \mathF _n$, all in the class $\CalH \CalT _{_{s}}$.
These
factors have at most countable fundamental group
(trivial in the case
$n < \infty$).
Moreover, in the case $n=\infty$ they can all be taken group von
Neumann algebras of the form $L(\mathF_\infty \ltimes \mathZ^2)$,
with $\mathF_\infty$ acting on $\mathZ^2$ as subgroups of
automorphisms via various embeddings $\mathF_\infty
\hookrightarrow \mathrm{SL}(2, \mathZ)$.
\end{corollary}

\noindent {\it Proof}. By the above comments and [Po01], all
actions $\alpha_t$ of $\mathF _n, 2\leq n \leq \infty,$ considered
in the proof of Corollary 5, give rise to factors $M_t =
L^\infty(\mathT ^2, \lambda) \rtimes_{\alpha_t} \mathF _n$ in the
class $\CalH \CalT _{_{s}}$. By Proposition~\ref{prop:4}, not only
the actions $\alpha_t$ are not stably orbit equivalent (modulo countable
sets), but even the factors $M_t$ are non stably isomorphic and
they all have countable fundamental group. Since for $n < \infty$
by [8.6.2$^\circ$, Po01] they have the same ($\neq 0, \infty$)
first $\ell^2_{_{HT}}$-Betti number, $\beta_1^{^{HT}}(M_t)=n-1,
\forall t$, the fundamental group $\CalF(M_t)$ is in this case
trivial, $\forall t$.

To see that in case $n=\infty$ the factors $M_t$ can be taken
group factors, write $\mathF_\infty$ as $\mathF_I$, the free group
with countable set of generators labeled by $I=(0,1)\cap \mathQ $,
as in Proposition 1. Then let $I_{t}=I\cap (0,t]$, for $t\in
(0,1],$ and denote by $\mathF_{I_t}$ the subgroup of $\mathF_I$
generated by generators with indices in $I_t$. Choose $\mathF_I\subset \mathrm{SL}(2, \mathZ)$ an embedding of $\mathF_I$ as a group of
automorphisms of $\mathZ^2$. By the remarks above again, the
factors $L(\mathF_{I_t} \ltimes \mathZ^2)$ are in the class $\CalH
\CalT _{_{s}}$ and by Proposition 4 they all have countable
fundamental group and are mutually non-stably isomorphic, modulo
countable sets. \hfill Q.E.D.

\begin{center} $\S $3.  \end{center}

The families of uncountably many non-equivalent free ergodic
p.m.p. actions of $\mathF_n$ in the previous section are all
obtained by approximating an initial equivalence relation $\CalS$
implemented by a relative property (T) action $\alpha$ of $\mathF_n$ from bellow by a one parameter family of sub-equivalence
relations $\CalS_t$ implemented  by actions $\alpha_t$ of the same
group $\mathF_n$: by the rigidity result ($4.5$ in \cite{Po01})
the relations $\CalS_t$ (and thus the actions $\alpha_t$) have the
relative property (T) when close enough to $\CalS$, and thus
Proposition 4 applies to obtain that modulo countable sets the
$\alpha_t$ are non-equivalent (even non stably orbit equivalent).

In this section we prove two more results in this vein. In both
cases, from an initial relative property (T) action of a group
$\Gamma$ we will produce new actions having this property, but of
groups which in general are different from the initial one. The
first result is based on a construction similar to the one in
Proposition ~\ref{prop:1}. The second one gives a general
hereditary property, showing that if an action $\sigma$ of a group
$\Gamma_0$ has the relative property (T) and $\Gamma_0 = \Gamma
\rtimes G$ with $G$ amenable, then $\sigma{|\Gamma}$ has this
property as well.

\begin{proposition}\label{prop:9} Let ${\alpha}$ be a free
p.m.p. action of ${\Gamma_{1}} * {\Gamma_{2}}$, with
${\Gamma_{1}}$ an infinite amenable group acting ergodically.
There exists a strictly increasing family of ergodic equivalence
relations ${\CalS }_{t}$, $t \in (0,1]$, implemented by free
ergodic p.m.p. actions $\alpha_t$ of the group $\Lambda=\Gamma_1 *
(*_{j\in \mathN } \Gamma_2)$, such that the associated Cartan
subalgebra inclusions define a strictly increasing family of
subfactors $A \subset M_t \subset M_1$, $0 < t \leq 1$, with $M_t
= A \rtimes_{\alpha_t} \Lambda$, $\forall t\in (0,1)$, $M_1= A
\rtimes_{\alpha} (\Gamma_1
* \Gamma_2)$, and $\overline{\cup_{t<s} M_t}^w= M_s$,  for each
$s\in (0,1]$. If in addition $\alpha$ has the relative property
$(T)$ (equivalently, if $A \subset M$ is rigid) then there exists
$c<1$ such that $\alpha_t$ has the relative property $(T)$
(equivalently $A \subset M_t$ is rigid) for all $t \in [c,1]$, and
thus the conclusions in $4.2^\circ$ hold true.
\end{proposition}

\smallskip
\noindent {\it Proof}. We prove a series of statements similar to
Proposition~\ref{prop:1}, keeping the hypotheses and the notation
$I,I_{t},G,G_{t},\Gamma_{1}, \Gamma_{2}, \alpha, \sigma_1$ from
there. For $t\in (0,1)$, define $\Lambda_t$ to be the kernel of
the natural epimorphism $G*\Gamma_{2}\to G/G_{t}\simeq
\oplus_{i\in I, i\not\in I_{t}} \mathZ /2 \mathZ $. Clearly,
$G*\Gamma_{2}$ is the strictly increasing union of the normal
subgroups $\Lambda_t$, $t\in (0,1)$. The restrictions
$\tilde{\sigma}_{t}$ of the action $\sigma_{1}=\tilde{\sigma}_{1}$
to the subgroups $\Lambda_{t}$ are ergodic on $G_t$ and define a
strictly increasing family of ergodic equivalence relations
$\tilde {\CalS }_{t}$, with the same properties as ${\CalS }_{t}$.

On the other hand, $\Lambda_t\simeq G_t* (*_{j\in \mathN }
\Gamma_2)$: the quotient of the Bass-Serre tree of $G*\Gamma_{2}$
by $\Lambda_{t}$ is a graph of groups made of a set of vertices of
degree 1, indexed by $g\in G/G_{t}$ with vertex groups $g^{-1}
\Gamma_{2}\  g$, all connected to one  ``central'' vertex with
vertex group $G_t$. By lemma~\ref{lem:2}, $\tilde{\sigma}_{t}$ is
OE to some ergodic free action $\tilde{\alpha}_{t}$ of $\Gamma_1 *
(*_{j\in \mathN } \Gamma_2)$.
  \hfill Q.E.D.

\vskip .05in

Observe that, with the notations in the above proof, by normality
of $\Lambda_{t}$, the outer automorphism group of $\CalR
_{\tilde{\alpha}_{t}, \Lambda}$ contains $G/G_{t} \simeq \oplus_{n
\in \mathN} \mathZ /2 \mathZ $ for all $t\in [c, 1)$. Notice also
that if we take $\Gamma_1\simeq \Gamma_2\simeq \mathZ $ then one
has $\Lambda_t\simeq G_t* \mathF _{\infty}$ and $\Gamma_1
* (*_{j\in \mathN } \Gamma_2)\simeq \mathF _{\infty}$.

\vskip .1in

\begin{proposition}\label{prop: 10} Let $\Gzero$ be a group of
the form $\Gzero= \Amen \ltimes \Gamma $, with $\Amen$ amenable.
Let $\sigma$ be a free p.m.p. action of $\Gzero$ on $(X, \mu)\Delta $.
Then $\sigma$ has the relative property $({\text{\rm T}})$ iff
$\sigma_{|\Gamma}$ has this property.
\end{proposition}

\begin{corollary} Let $\Gzero$ be a group of
the form $\Gzero=\Amen \ltimes \Gamma$, with $\Amen$ amenable.
Assume $\Gzero$ acts by outer automorphisms on an abelian group
$\Delta $. Then the pair $(\Delta  \ltimes \Gzero, \Delta )$ has the relative
property (T) iff $(\Delta  \ltimes \Gamma, \Delta )$ has this property.

\end{corollary}

\noindent {\it Proof}. By (5.1 in [Po01]), it is sufficient to
prove Proposition 10. Let $A=L^\infty(X, \mu)$, $\Mzero= A
\rtimes_\sigma \Gzero$, $M = A \rtimes_\sigma \Gamma$. If
$A\subset M$ is rigid then $A \subset \Mzero$ is clearly rigid.

Conversely, assume $A \subset \Mzero$ is rigid and  let
$\Fbig=\Fbig(\varepsilon/16)\subset \Mzero$,
$\deltabig=\deltabig(\varepsilon/16)> 0$ be so that if $\phibig$
is a completely positive map on $\Mzero$ with $\tau\circ \phibig
\leq \tau$, $\phibig(1) \leq 1$ and $\|\phibig(x)-x\|_2 \leq
\deltabig$ for any $x \in \Fbig$ then $\|\phibig(a)-a\|_2 \leq
\varepsilon/16, \forall a\in A$, $\|a\|\leq 1$.

By the amenability of $\Amen$, it follows that there exists a
finite subset $K\subset \Amen$ such that
\begin{eqnarray}
\|\ \frac{1}{|K|} \Sigma_{g,h\in K} u_gE_M(u_g^*xu_h)u_h^* - x
\|_2 \leq \deltabig/2, \forall x\in \Fbig
\end{eqnarray}
By using again the amenability of $\Amen$, by the Ornstein-Weiss
Theorem [OW80] it follows that for any $\deltabig > 0$ there
exists a finite group $\Amen_{0}=\oplus_j (\mathZ /2\mathZ )_j$
and unitary elements $\{v_g\}_{g\in \Amen_{0}}$ in the normalizer
of $A$ in $A \rtimes \Amen \subset \Mzero$ which implement a free
action of $\Amen_{0}$ on $(X, \mu)$ such that each $u_g, g \in K$,
is $\varepsilon/64$-close to an element $u'_g$ in the normalizer
of $A$ in $A \vee \{v_g\}_{g\in \Amen_{0}}\simeq A \rtimes
\Amen_{0}$.

Thus, if we denote by $A_{0}=A\cap \{v_g\}'_g$ the fixed point
algebra of the action of $\Amen_{0}$ on $A$ implemented by the
unitaries $v_g$, then $A$ can be decomposed as $A=A_{0} \otimes
A_1$ with $A_1$ generated by unitary elements $w_h, h \in
\hat{\Amen}_{0}$, satisfying $v_gw_hv_g^* = h(g) w_h$, $\forall g
\in \Amen_{0}, h \in \hat{\Amen}_{0}$.

There exists $\delta > 0$ such that if $\phi$ is a completely
positive map on $M$ with $\tau \circ \phi \leq \tau$, $\phi(1)
\leq 1$ and $\|\phi(w_h)-w_h\|_2 \leq \delta, \forall h\in
\hat{\Amen}_0,$  then $\|\phi(uw_h)- \phi(u)w_h\|_2 \leq
\varepsilon/16$, $\forall u\in \CalU (A), h\in \hat{\Amen}_{0}$
(see e.g., Corollary 1.1.2 in [Po01]). Thus, for $u \in \CalU
(A_{0}), h \in \hat{\Amen}_{0}$ we have
\begin{eqnarray}
\|\phi(uw_h)-uw_h\|_2 \leq \|\phi(uw_h) -\phi(u)w_h\|_2 +
\|\phi(u) -u\|_2 \leq \varepsilon/16 + \|\phi(u) -u\|_2
\end{eqnarray}
Moreover, $\delta$ can be chosen small enough so that we also have
\begin{eqnarray}
\| \phi(E_M(u_g^*xu_{g'})) - E_M(u_g^*xu_{g'})\|_2 \leq
(2|K|)^{-1}\deltabig, \forall g,g' \in K
\end{eqnarray}

We set $F(\varepsilon)= \{E_M(u_g^*xu_{g'}) \mid x \in
\Fbig(\varepsilon/16), g, g' \in K\}\cup \{w_h \mid h \in
\hat{\Amen}_{0}\}$, $\delta(\varepsilon)=\text{\rm min} \{\delta,
\deltabig, \varepsilon/4\}$ and will show that if $\phi$ is a
completely positive map on $M$ such that $\tau \circ \phi \leq
\tau$, $\phi(1)\leq 1$ and $\|\phi(x)-x\|_2 \leq
\delta(\varepsilon), \forall x\in F(\varepsilon)$ then
$\|\phi(a)-a\|_2 \leq \varepsilon, \forall a\in A$, $\|a\|\leq 1$.

To see this, note first that if we put $\phibig(x)
\overset{\mathrm{def}}{=}|K|^{-1}\Sigma_{g,g' \in K}
u_g\phi(E_M(u_g^*xu_{g'}))u_{g'}^*$ for $x\in \Mzero$ then
$\phibig$ is a completely positive map on $\Mzero$ with $\tau
\circ \phibig \leq \tau, \phibig(1) \leq 1$. Moreover,  by $(1)$
and $(3)$, for each $x \in \Fbig(\varepsilon/16)$ we have the
inequalities:
$$
\|\phibig(x)-x\|_2 \leq |K|^{-1} \Sigma_{g,g'\in K} \|
u_g\phi(E_M(u_g^*xu_{g'}))u_{g'}^*-
u_gE_M(u_g^*xu_{g'})u_{g'}^*\|_2
$$
$$
+ \||K|^{-1}\Sigma_{g,g'\in K} u_gE_M(u_g^*xu_{g'})u_{g'}^* -
x\|_2 \leq \deltabig/2 + \deltabig/2 = \deltabig.
$$

Thus, $\|\phibig(a)-a\|_2 \leq \varepsilon/16$, $\forall a \in A$,
$\|a\|\leq 1$. Since for $a \in A$ we have
$$
\phibig(a)=|K|^{-1} \Sigma_{g\in K} u_g\phi(E_M(u_g^*au_g))u_g^*=
|K|^{-1} \Sigma_{g\in K} u_g \phi(\sigma_g^{-1}(a))u_g^*
$$
and $\|u_g -u'_g\|_2 \leq \varepsilon/64$, this implies that for
$a_{0} \in A_{0}, \|a_{0}\|\leq 1$ we have:
$$
\| |K|^{-1} \Sigma_g u'_g\phi(a_{0}){u'}_g^* - a_{0}\|_2 \leq
\|\phibig(a_{0})-a_{0}\|_2 + 4\|K\|^{-1}\Sigma_{g\in K}
\|u_g-u'_g\|_2
$$
$$
\leq \varepsilon/16 + \varepsilon/16 = \varepsilon/8.
$$

By convexity, since
$\|u'_g\phi(a_{0}){u'}_g^*\|_2=\|\phi(a_{0})\|_2 \leq \|a_{0}\|$,
it follows that for ``most'' $g \in K$ (one in fact is sufficient
!) we have
$$
\|\phi(a_{0})-a_{0}\|_2=\|u'_g\phi(a_{0}){u'}_g^*-a_{0}\|_2 \leq
\varepsilon/8, \forall a_{0} \in A_{0}, \|a_{0}\|\leq 1.
$$

Thus, with the notations in (1.1.2 of [Po01]), by using part
$4^\circ$ of (Lemma 1.1.2 in [Po01]), in the Hilbert $M$-bimodule
$\CalH _\phi$ we have $\|[u_{0}, \xi_\phi]\|_2 \leq
2\|\phi(u_{0})-u_{0}\|_2\leq \varepsilon/4$, $\forall u_{0} \in
\CalU (A_{0})$ and $\|[w_h, \xi_\phi]\|_2 \leq \varepsilon/4$,
$\forall h \in \hat{\Amen}_{0}$, as well. This implies $\xi_\phi$
is $\varepsilon/2$-close to a vector $\xi \in \CalH _\phi$ that
commutes with both $\CalU (A_{0})$ and the group $\{w_h\}_h$.

Since $\CalU (A_{0})$ and $\{w_h\}_h$ generate $A$, it follows
that $\xi$ commutes with all the elements in $A$. Thus, $\|[a,
\xi_\phi]\|_2 \leq \varepsilon/2$, $\forall a \in A$, $\|a\| \leq
1$. By (part $2^\circ$ of Lemma 1.1.3 in [Po01]), this implies
$\|\phi(a)-a\|_2 \leq \varepsilon, \forall a\in A, \|a\|\leq 1$.
\hfill Q.E.D.

\vskip .05in \noindent {\bf Final Remarks}. $1^\circ$. Note that
Propositions 9 and 4, taken together, provide another construction
of uncountable families of non stably orbit equivalent actions of $\mathF_\infty$. As concerning orbit equivalence, observe that it is in
fact enough, by \cite[Prop. II.6]{Ga00} and \cite{Hj02-a}, to get one
free ergodic p.m.p. action of
$\mathF _{\infty}$ with countable fundamental group to get
uncountably many non OE actions by restriction to Borel subsets of
various measures. However, these actions are stably orbit equivalent.

$2^\circ$. Note that in Proposition 10 we have actually proved a
slightly more general statement, namely: Let $M$ be a finite von
Neumann algebra with $A \subset M$ a relatively rigid Cartan
subalgebra. Let $A \subset N \subset M$ be an intermediate von
Neumann subalgebra such that $M=N \rtimes_\sigma H$ for some free
action of an amenable group $H$ on $N$, with the property that
$\sigma_h(A)=A, \forall h\in H$. Then $A$ is relatively rigid in
$N$ as well. Related to this statement, it would be interesting to
investigate the following question: {\it Assuming $B \subset N
\subset M$ is an inclusion of finite von Neumann algebras such
that $M$ is amenable relative to $B$, in the sense of Definition
3.5.2$^\circ$ of} [Po01], {\it and $B \subset M$ is rigid, does $B
\subset N$ follow rigid as well ?}

$3^\circ$. Let $\CalS$ be the equivalence relation implemented by
the action of $\mathrm{SL}(2, \mathZ)$ on $(\mathbb T^2, \lambda)$. By
\cite{Bu91} and \cite{Po01}, any sub-equivalence relation
$\CalS_0$ of $\CalS$ implemented by a nonamenable subgroup
$\Gamma_0 \subset \mathrm{SL}(2, \mathZ)$ has the relative property (T). By
the rigidity result (5.4 in \cite {Po01}), Propositions 1, 4, 9
and 10 the class of subequivalence relations $\CalS_0 \subset
\CalS$ for which this still holds true gets considerably larger.
The following problem thus seems quite natural: {\it Do all
nonamenable ergodic subequivalence relations $\CalS_0$ of $\CalS$
have the relative property} (T) ?

\bigskip

\medskip

\noindent {\bf Acknowledgments}. We are grateful to \'E. Ghys for
helpful indications concerning Lemma~\ref{lem:6} and to B. Bekka
for pointing out to us the results on relative property (T) in
[Bu91].

\bigskip
\noindent \textsc{D.~G.: UMPA, UMR CNRS 5669, ENS-Lyon, F-69364
Lyon
Cedex 7}

\noindent
\textsc{S.~P.: Math.Dept., UCLA, LA, CA 90095-155505, USA}

\noindent \texttt{gaboriau@umpa.ens-lyon.fr, \ popa@math.ucla.edu}

\end{document}